\documentclass[12pt]{article}
\usepackage{amsfonts}

\usepackage{latexsym}
\usepackage{amssymb}
\usepackage{amsmath}
\usepackage{amsthm}
\usepackage[mathscr]{eucal}
\usepackage{graphicx}
\usepackage{hyperref}
\usepackage{caption}
\usepackage{subcaption}

\newtheorem{Lemma}{Lemma}

\newtheorem{Theorem}[Lemma]{Theorem}

\newtheorem{Definition}{Definition}

\renewcommand{\qed}{\hfill{\ \ \rule{2mm}{2mm}} \vspace{0.2in}}

\begin{document}

\title{Strong Identification Codes for Graphs}
\author{ \textbf{Ghurumuruhan Ganesan}
\thanks{E-Mail: \texttt{gganesan82@gmail.com}} \\
\ \\
Institute of Mathematical Sciences, HBNI, Chennai}
\date{}
\maketitle


\begin{abstract}
For any graph~\(G,\) a set of vertices~\({\cal V}\) is said to be dominating if every vertex of~\(G\) contains at least one node of~\(G\) and separating if each vertex~\(v\) contains a unique neighbour~\(u_v \in {\cal V}\) that is adjacent to no other vertex of~\(G.\) If~\({\cal V}\) is both dominating and separating, then~\({\cal V}\) is defined to be an identification code. In this paper, we study strong identification codes with an index~\(r,\) by imposing the constraint that each vertex of~\(G\) contains at least~\(r\) unique neighbours in~\({\cal V}.\) We use the probabilistic method to study both the minimum size of strong identification codes and the existence of graphs that allow an identification code with a given index.


\vspace{0.1in} \noindent \textbf{Key words:} Strong identification codes, strong neighbourhood graphs, probabilistic method.

\vspace{0.1in} \noindent \textbf{AMS 2000 Subject Classification:} Primary: 05C62;
\end{abstract}

\bigskip

\renewcommand{\theequation}{\thesection.\arabic{equation}}
\setcounter{equation}{0}
\section{Introduction} \label{intro}

An identification code of a graph~\(G\) is a set of vertices of~\(G\) that form both a dominating and a separating set (Karpovsky, Chakrabarty and Levitin~\cite{karpo}). Bounds for~\(\gamma(G),\) the minimum size of an identification code for~\(G,\) have been evaluated for different classes of graphs: Gravier and Moncel~\cite{grav} describe graphs whose minimum size identification code contains all vertices apart from exactly one vertex. Foucaud, Klasing, Kosowski and Raspaud~\cite{flor2} studied constructions of identification codes for triangle-free graphs, with size at most a constant fraction of the number of vertices. Later Foucaud and Perarnau~\cite{flor1} obtained bounds for~\(\gamma(G)\) in terms of the minimum vertex degree~\(\delta\) for graphs~\(G\) whose girth is at least five and also study the size of the identification codes for random regular graphs. Recently,~\cite{balb} have investigated the use of graph spectra as a tool to characterize the identification codes.


In this paper, we study strong identification codes of a graph~\(G\) with an index~\(r \geq 1,\) satisfying the additional constraint that each vertex contains at least~\(r\) unique neighbours in the code. (For formal definitions, we refer to Section~\ref{sec_str_code}).  We use the probabilistic method to estimate the minimum size of an identification code and also to determine existence of graphs that allow strong identification codes with a given index.

The paper is organized as follows. In Section~\ref{sec_str_code}, we state and prove our main result, Theorem~\ref{thm_iden}, regarding the minimum size of a strong identification code with index~\(r.\) Next in Section~\ref{sec_str_neigh}, we state and prove Theorem~\ref{prop1} regarding existence of graphs that allow strong idenfication codes with a given index.



\renewcommand{\theequation}{\thesection.\arabic{equation}}
\setcounter{equation}{0}
\section{Strong Identification Codes} \label{sec_str_code}
Let~\(G = (V,E)\) be a graph with vertex set~\(V = \{1,2,\ldots,n\}\) and edge set~\(E.\) For a vertex~\(v \in G,\) we define~\(N(v)\) to be the set of all vertices adjacent to~\(v,\) not including~\(v\) and let~\(N[v] = N(v) \cup \{v\}.\) We also denote the elements of~\(N(v)\) as \emph{neighbours} of~\(v.\) Letting~\(\#A\) denote the number of elements in the set~\(A,\) we have the following definition.
\begin{Definition}
A set of vertices~\({\cal C}\) is said to be an \emph{identification code} for~\(G\) with \emph{index}~\(r\) if for any two distinct vertices~\(v,u \in G\) we have
\begin{equation}\label{ghj}
\#\left(\left(N[v] \setminus N[u]\right) \bigcap {\cal C}\right) \geq r.
\end{equation}
\end{Definition}
For~\(r=1,\) the above definition reduces to the concept of identification codes  as studied in Gravier and Moncel~\cite{grav} and Foucaud and Perarnau~\cite{flor1}.

Clearly, one necessary condition for~\(G\) to have an identification code with index~\(r\) is that the neighbourhood of any vertex contains at least~\(r\) distinct vertices not present in any other vertex, i.e., for any two distinct vertices~\(u,v \in G,\) we must have that
\begin{equation}\label{iden_eq}
\#\left(N[v] \setminus N[u] \right) \geq r.
\end{equation}
If~(\ref{iden_eq}) holds, we say that~\(G\) has the~\(r-\)strong neighbourhood property. In the next Section~\ref{sec_str_neigh}, we use the probabilistic method to establish the existence of graphs satisfying~(\ref{iden_eq}).

If the~\(r-\)strong neighbourhood property holds, then the whole vertex set of~\(G\) is itself an identification code with index~\(r.\) For graphs with stronger neighbourhoods, we have the following result regarding the minimum size of an identification code with a given index.
\begin{Theorem}\label{thm_iden} Let~\(G\) be a graph on~\(n\) vertices satisfying the~\((r+d+1)-\)strong neighbourhood property for some integer~\(d \geq 1.\) If~\(\theta_n = \theta_n(G,r,d)\) is the minimum size of an identification code of~\(G\) with index~\(r,\) then
\begin{equation}\label{theta_bounds}
\frac{1}{\Delta+1} \leq \frac{\theta_n}{n} \leq 1-\frac{c(d,r)}{(\Delta+1)^{(r+2)/d}},
\end{equation}
where~\(\Delta \geq 2\) is the maximum degree of a vertex in~\(G\) and\\\(c(d,r) := d(d+1)^{-1-1/d}(2r)^{-1/d}.\)
\end{Theorem}
Thus if~\(G\) has~\((r+d+1)-\)strong neighbourhood property for some integer~\(d \geq 1,\) then there exists a \emph{nontrivial} fraction of vertices of~\(G\) that act as an identification code with index~\(r.\)


\subsection*{Proof of Theorem~\ref{thm_iden}}
The lower bound in~(\ref{theta_bounds}) is obtained as follows. If~\({\cal C}\) is an identification code with index~\(r,\) then~\({\cal C}\) is also a dominating set in the sense that~\(\bigcup_{v \in {\cal C}} N[v] = \{1,2,\ldots,n\}.\) Since~\(\#N[v] \leq \Delta+1,\) it follows that~\(\#{\cal C} \geq \frac{n}{\Delta+1}.\)

To prove the upper bound in~(\ref{theta_bounds}), we use the probabilistic method similar to Theorem~1.2.2, pp. 4, Alon and Spencer~\cite{alon}, that concerns the minimum size of a dominating set. We select each vertex of~\(G\) independently with probability~\(q\) to be determined later and let~\(Z\) be the random set of vertices selected.
Say that a vertex~\(v \in G\) is \emph{bad} if one of the following two conditions hold: either~\((a)\) \(\#\left(N[v] \bigcap Z \right) \leq r-1\) or~\((b)\)~\(\#\left(N[v] \bigcap Z\right) \geq r\) but there exists a vertex~\(w\) satisfying
\begin{equation}\label{gen_h}
\#\left((N[v]\setminus N[w]) \bigcap Z\right) \leq r-1.
\end{equation}
A vertex of~\(G\) which is not bad is said to be \emph{good}.

If~\(Y_b\) is the set of all bad vertices and~\(Z_b = \bigcup_{v \in Y_b} N[v],\) then the set~\(Y := Z \cup Z_b\) is an identification code for~\(G\) with index~\(r:\) Indeed, if~\(v\) is a good vertex, then for any vertex~\(u \neq v\) we have
\[\#\left((N[v] \setminus N[u]) \cap Y\right) \geq \#\left((N[v] \setminus N[u])\cap Z\right) \geq r.\]
On the other hand if~\(v\) is a bad vertex, then for any vertex~\(u \neq v,\) we have by the~\((r+d+1)-\)strong neighbourhood property that
\begin{eqnarray}
\#\left((N[v] \setminus N[u]) \cap Y\right) &\geq& \#\left((N[v] \setminus N[u]) \cap Z_b\right)  \nonumber\\
&=& \#\left(N[v] \setminus N[u]\right) \nonumber\\
&\geq& r+d+1. \nonumber
\end{eqnarray}


We now estimate the size of~\(Z_b\) by bounding the probability that a vertex~\(v \in G\) is bad. If~\(\delta(v) := \#N[v],\) then the probability that case~\((a)\) occurs is
\begin{equation}\label{case_a}
f_1(q,v) := \sum_{l=0}^{r-1} {\delta(v) \choose l} q^{l}(1-q)^{\delta(v)-l}.
\end{equation}
Now if possibility~\((b)\) occurs and some vertex~\(w\) satisfies~(\ref{gen_h}), then either~\(w\) is a neighbour of~\(v\) or~\(w\) is at distance two from~\(v\) and shares a common neighbour with~\(v.\) Letting~\(\delta(v,w) := \#(N[v] \setminus N[w]),\) the probability that case~\((b)\) occurs is then bounded above by
\begin{equation}
f_2(q,v) := \max_{w} \sum_{l=0}^{r-1} {\delta(v,w) \choose l} q^{l}(1-q)^{\delta_1(v,w)-l},
\label{case_b_one}
\end{equation}
where the maximum is taken over all vertices~\(w\) at a distance of at most two from~\(v.\) Since the maximum degree of any vertex is at most~\(\Delta,\) there are at most~\(\Delta + \Delta(\Delta-1) = \Delta^2\) neighbours at a distance of at most two from~\(v.\) Consequently, the probability that vertex~\(v\) is bad is upper bounded by~\(f_1(q,v) + \Delta^2 f_2(q,v)\) and so the expected number of vertices in~\(Y_b\) is at most~\(n( f_1(q,v) + \Delta^2 f_2(q,v)).\) This in turn implies that the expected number of vertices in~\(Z_b =\bigcup_{v \in Y_b}N[v]\) is at most~\(n(\Delta f_1(q,v) + \Delta^3 f_2(q,v)).\)

Summarizing, we get from the discussion in the previous paragraph that the expected number of vertices in~\(Z \cup Z_b\) is at most
\begin{equation}\label{bad_exp}
n\left(q + \Delta \max_{v} f_1(q,v) + \Delta^3 \max_{v} f_2(q,v)\right).
\end{equation}
To find suitable bounds for~\(f_1(q,v)\) and~\(f_2(q,v),\) we first let~\[\delta := \min_{v} \#N(v) \geq r+d\] be the minimum degree of a vertex in~\(G.\) From~(\ref{case_a}) we have that
\[f_1(q,v) \leq (1-q)^{\delta-r+1} \sum_{l=0}^{r-1} {\delta(v) \choose l} q^{l} \leq (1-q)^{d+1} \sum_{l=0}^{r-1} {\delta(v) \choose l} q^{l}\] and using~\({\delta(v) \choose l} \leq \delta^{l}(v) \leq (\Delta+1)^{l}\) and~\(q \leq 1,\) we get that
\begin{equation}\label{f_one_up}
f_1(q,v)\leq (1-q)^{d+1} \sum_{l=0}^{r-1}(\Delta+1)^{l} \leq r(\Delta+1)^{r-1}  (1-q)^{d+1}.
\end{equation}

For estimating~(\ref{case_b_one}), we argue similarly and use~\(q \leq 1\) and the fact that~\[\delta(v,w)  = \#(N[v] \setminus N[w]) \geq r+d+1\] by the~\((r+d+1)-\)strong neighbourhood property, to get that
\begin{equation}\label{f_2_up}
f_2(q,v) \leq (1-q)^{d+2} \sum_{l=0}^{r-1} {\delta(v,w) \choose l} \leq r(\Delta+1)^{r-1} (1-q)^{d+2}.
\end{equation}
Plugging the estimates~(\ref{f_one_up}) and~(\ref{f_2_up}) into~(\ref{bad_exp}), we get that the expected number of vertices in the identification code~\(Z \cup Z_b\) is at most~
\begin{eqnarray}
\mathbb{E} \left(\#(Z \cup Z_b)\right) &\leq& n q+ r(\Delta+1)^{r} (1-q)^{d+1} + r(\Delta+1)^{r+2} (1-q)^{d+2} \nonumber\\
&\leq& n(q+ 2r(\Delta+1)^{r+2}(1-q)^{d+1}) \nonumber\\
&=:& n \Gamma(q).
\end{eqnarray}

The solution~\(q_0\) to the equation~\(\Gamma'(q) = 0\) satisfies
\[1-2r(\Delta+1)^{r+2}(d+1)(1-q_0)^{d} =0\] and so~\[\Gamma(q_0) =q_0 + \frac{1-q_0}{d+1} =  1-\frac{d}{d+1}  \left(\frac{1}{2r(\Delta+1)^{r+2}(d+1)}\right)^{\frac{1}{d}}.\]
Consequently, there exists at least one identification code with index~\(r\) and size at most~\(n \Gamma(q_0).\)~\(\qed\)

\renewcommand{\theequation}{\thesection.\arabic{equation}}
\setcounter{equation}{0}
\section{Strong Neighbourhood Graphs} \label{sec_str_neigh}
In this section, we use the probabilistic method to obtain~\(r-\)strong neighbourhood graphs for any integer~\(r \geq 1.\)

The cycle~\(C_n\) on~\(n \geq 4\) vertices satisfies the~\(1-\)strong neighbourhood property. The following result establishes the existence of bounded degree graphs with~\(w-\)strong neighbourhoods, for any integer~\(w\geq 2.\)
\begin{Theorem}\label{prop1} Let~\(w \geq 2\) be fixed. There exists an absolute constant~\(C > 0\) not depending on~\(w\) such that if~\(n \geq \max(C,160(w+1)^2+1) := M(w),\) then there exists a connected graph~\(G\) on~\(n\) vertices satisfying the~\(w-\)strong neighbourhood property and having maximum vertex degree~\[\Delta(G) \leq \max(10\log(2M(w)),8(w+1))+1.\]
\end{Theorem}
Thus for any~\(w \geq 2\) and all~\(n\) large, there are bounded degree graphs on~\(n\) vertices satisfying the~\(w-\)strong neighbourhood property.

We use the probabilistic method to prove Theorem~\ref{prop1}, beginning with the following auxiliary Lemma.
\begin{Lemma}\label{deg_lem}
Let~\(y \geq 3\) be fixed. There exists an absolute constant~\(C > 0\) not depending on~\(y\) such that if~\(n \geq \max(C,160y^2+1),\) then there exists a connected graph~\(G\) satisfying the~\((y-1)-\)strong neighbourhood property and having maximum vertex degree of at most~\(\max(32\log{n},8y).\)
\end{Lemma}

To prove Lemma~\ref{deg_lem}, we use the following standard deviation estimate. Let~\(\{X_j\}_{1 \leq j \leq L}\) be independent Bernoulli random variables with~\(\mathbb{P}(X_j = 1) = 1-\mathbb{P}(X_j = 0) > 0.\) If~\(T_L = \sum_{j=1}^{L} X_j,\theta_L = \mathbb{E}T_L\) and~\(0 < \epsilon \leq \frac{1}{2},\) then
\begin{equation}\label{conc_est_f}
\mathbb{P}\left(\left|T_L - \theta_L\right| \geq \theta_L \epsilon \right) \leq 2\exp\left(-\frac{\epsilon^2}{4}\theta_L\right).
\end{equation}
For a proof of~(\ref{conc_est_f}), we refer to Corollary A.1.14, pp. 312 of Alon and Spencer (2008).

\emph{Proof of Lemma~\ref{deg_lem}}: We use the probabilistic method. Let~\(H = G(n,p)\) be the random graph obtained by setting each edge of the complete graph~\(K_n\)
on~\(n\) vertices to be independently present with probability~\(p,\) where
\begin{equation}\label{p_def}
p = \frac{\max(16\log{n},4y)}{n-1}
\end{equation}
and absent otherwise.

Letting~\(N(1)\) denote the neighbours of the vertex~\(1,\) we have that\\\(\mathbb{E}\#N(1) = (n-1)p\) and so using~(\ref{conc_est_f}) with~\(\epsilon = \frac{1}{2},\) we get
\begin{equation}\label{n_one_est}
\mathbb{P}\left(\frac{(n-1)p}{2} \leq \#N(1) \leq \frac{3(n-1)p}{2}\right) \geq 1-\exp\left(-\frac{(n-1)p}{4}\right) \geq 1-\frac{1}{n^4},
\end{equation}
by~(\ref{p_def}). Letting
\[A_{tot} := \bigcap_{1 \leq i \leq n} \left\{\frac{(n-1)p}{2} \leq \#N(1) \leq \frac{3(n-1)p}{2}\right\},\] we get by the union bound and~(\ref{n_one_est}) that
\begin{equation}\label{a_tot_est}
\mathbb{P}(A_{tot}) \geq 1-\frac{1}{n^3}.
\end{equation}
Thus, with high probability, i.e., with probability converging to one as~\(n~\rightarrow~\infty,\) each vertex in the random graph~\(H\) contains at most~\(2(n-1)p = \max(32\log{n}, 8y)\) neighbours.

Next, we estimate the number of common neighbours between the two vertices~\(i\)
and~\(j\) as follows. Letting~\(T_{ij} := \#(N(i) \cap N(j))\) and using the Chernoff bound, we have for~\(s,t >0\)
that
\begin{equation}\label{eq_one}
\mathbb{P}\left(T_{ij} > t\right) \leq \left(\mathbb{E}e^{sT_{ij}}\right)^{n-2}e^{-st}.
\end{equation}
Each vertex~\(z \neq i,j\) is adjacent to both~\(i\) and~\(j\) with probability~\(p^2,\) independent of the other vertices and so
\begin{equation}\label{eq_two}
1 < \mathbb{E}e^{sT_{ij}} = 1 + (e^{s}-1)p^2 \leq \exp\left((e^{s}-1) p^2\right).
\end{equation}
Substituting~(\ref{eq_two}) into~(\ref{eq_one}) and using~\((\mathbb{E}e^{sT_{ij}})^{n-2} \leq (\mathbb{E}e^{sT_{ij}})^{n-1},\) we get
\begin{equation}\label{eq_3}
\mathbb{P}\left(T_{ij} > t\right) \leq \exp\left((e^{s}-1) (n-1)p^2\right)e^{-st}.
\end{equation}

For~\(s=2,\) we have
\begin{eqnarray}
(e^{s}-1)(n-1)p^2 &=& (e^{2}-1)\frac{\max(256(\log{n})^2,16y^2)}{n-1} \nonumber\\
&\leq& \frac{1}{n-1}\max(2560(\log{n})^2,160y^2), \nonumber
\end{eqnarray}
which is at most one, provided that~\(n \geq 1+\max(160y^2,n_0),\) where~\(n_0\) is the smallest integer~\(x\) such that~\(x > 2560(\log{x})^2.\)
For all~\(n \geq 1+\max(160y^2,n_0),\) we therefore have from~(\ref{eq_3}) that
\begin{equation}\label{eq_32}
\mathbb{P}\left(T_{ij} > \frac{(n-1)p}{4}\right) \leq e\exp\left(-\;\frac{(n-1)p}{2}\right) \leq \frac{e}{n^{8}},
\end{equation}
since~\((n-1)p \geq 16\log{n}\) by~(\ref{p_def}). Thus,
\[\mathbb{P}\left(\#\left(N(i) \cap N(j)\right) \leq \frac{(n-1)p}{4}\right)  \geq 1-\frac{e}{n^{8}}\]
and, using the union bound, we have
\begin{equation}\label{hu}
\mathbb{P}\left(\bigcap_{i \neq j} \left\{\#\left(N(i) \cap N(j)\right) \leq \frac{(n-1)p}{4}\right\}\right)  \geq 1-\frac{e}{n^{6}}.
\end{equation}

From~(\ref{a_tot_est}) and~(\ref{hu}), we get that with high probability, the maximum degree of the random graph~\(H\) is at most~\(2(n-1)p = \max(32\log{n},8y).\)
In addition, for any two vertices~\(v \neq u,\) the corresponding neighbourhood sets satisfy
\[\#\left(N[v] \setminus N[u]\right) \geq \frac{(n-1)p}{4}-1 \geq y-1.\]
Consequently, with high probability, the graph~\(H\) satisfies the~\((y-1)-\)strong neighbourhood property
and because~\(p\geq \frac{16\log{n}}{n}\) (see~(\ref{p_def})), the graph~\(H\) is also connected with high probability (pp. 164--165, Bollob\'as~\cite{boll}). Summarizing, there exists an absolute constant~\(C > 0\) not depending on~\(y\) such that if\\\(n \geq \max(C,160y^2~+~1),\) then with high probability~\(H\) is connected,
satisfies the~\((y-1)-\)strong neighbourhood property and has a maximum vertex degree of at most~\(\max(32\log{n},8y).\)~\(\qed\)

\emph{Proof of Theorem~\ref{prop1}}: Let~\(w \geq 2\) and~\(M(w)\) be as in the statement of Theorem~\ref{prop1}. For~\(n \geq M(w),\) split~\(\{1,2,\ldots,n\}\)
into~\(T\) sets~\(V_1,\ldots,V_T\) such that~\(V_i, 1 \leq i \leq T-1,\) has size exactly~\(M(w)\) and~\(V_T\) has size at least~\(M(w)\)
and at most~\(2M(w).\) Using Lemma~\ref{deg_lem} with~\(w = y-1,\) we know that there exists a connected graph~\(G_i\) with vertex set~\(V_i\) satisfying the~\(w-\)strong neighbourhood property and having a maximum vertex degree of at most\\\(\Delta_0 = \max(32\log(2M(w)),8(w+1)).\) We pick one vertex~\(u_i \in V_i\) and connect~\(u_i\) to~\(u_{i+1}\) for~\(1 \leq i \leq T-1.\) The resulting graph~\(G\) is connected, satisfies the~\(w-\)strong neighbourhood property and has a maximum vertex degree of at most~\(\Delta_0+1.\)~\(\qed\)



%





\subsection*{Acknowledgement}
I thank Professors Rahul Roy, C. R. Subramanian and the referees for crucial comments that led to an improvement of the paper. I also thank IMSc for my fellowships.

\bibliographystyle{plain}

\end{document}